\documentclass[twocolumn,aps,secnumarabic,nofootinbib,floatfix]{revtex4}
\usepackage{amssymb,amsmath2000,amsthm2000,times,mathptmx,pstricks}
\newtheorem{theorem}{Theorem}[section]
\newtheorem{lemma}[theorem]{Lemma}
\newtheorem{corollary}[theorem]{Corollary}

\newcommand{\C}{\mathbb{C}}
\newcommand{\Z}{\mathbb{Z}}

\newcommand{\SO}{\mathrm{SO}}
\newcommand{\Sp}{\mathrm{Sp}}
\newcommand{\SL}{\mathrm{SL}}
\renewcommand{\tensor}{\otimes}
\newcommand{\fig}[1]{{Figure~\ref{#1}}}

\newcommand{\lzef}{\pspolygon[fillstyle=solid,fillcolor=white]
    (0,0)(.866,.5)(1.732,0)(.866,-.5)}
\newcommand{\lzsef}{\pspolygon[fillstyle=solid,fillcolor=lightgray]
    (0,0)(0,-1)(.866,-1.5)(.866,-.5)}
\newcommand{\lzswf}{\pspolygon[fillstyle=solid,fillcolor=gray]
    (0,0)(0,-1)(-.866,-1.5)(-.866,-.5)}
\psset{unit=.333in,linewidth=.3pt,dash=2.5pt 2.5pt}

\begin{document}

\title{Self-complementary plane partitions by Proctor's minuscule method}
\author{Greg Kuperberg}
\thanks{Supported by NSF Postdoctoral Fellowship DMS \#9107908}
\affiliation{Department of Mathematics,
    University of Chicago, Chicago, IL 95616}
\email[Current email: ]{greg@math.ucdavis.edu}

\begin{abstract}
A method of Proctor \cite{Proctor:bruhat} realizes the set of arbitrary plane
partitions in a box and the set of symmetric plane partitions as bases of
linear representations of Lie groups. We extend this method by realizing
transposition and complementation of plane partitions as natural linear
transformations of the representations, thereby enumerating symmetric plane
partitions, self-complementary plane partitions, and transpose-complement
plane partitions in a new way.
\end{abstract}

\maketitle

While investigating minuscule representations of Lie algebras, Robert Proctor
\cite{Proctor:bruhat} developed one of the simplest methods for enumerating
symmetry classes of plane partitions. The method yields an enumeration of
arbitrary plane partitions, using the representation theory of $\SL(n)$, and
symmetric plane partitions, using the (projective) representation theory of
$\SO(2n+1)$.  (Proctor also provided solutions to two previously open symmetry
classes, namely  TCPP's and SSCPP's in the notation of Stembridge, in a related
representation-theoretic investigation.)

The author and John Stembridge have independently discovered a way to extend
this method to enumerate self-complementary plane partitions (first enumerated
by Stanley \cite{Stanley:symmetries}).  In this paper, we will present a brief
proof of this enumeration; Stembridge \cite{Stembridge:minuscule} gives a
longer exposition and another generalization of the new proof.  For both of us,
the proof was inspired by, and partially explains, Stembridge's $q = -1$
phenomenon \cite{Stembridge:q-1}, which in this case is the observation that
the number of self-complementary plane partitions is the number of plane
partitions with an even number of cubes minus the number with an odd number of
cubes.

In addition, the author has found a new proof of the enumeration of symmetric
plane partitions which extends to an enumeration of TCPP's,  while Stembridge's
analysis also enumerates SSCPP's and self-complementary chains of order ideals
in other minuscule posets.  Unfortunately, symmetry classes involving cyclic
symmetry have no known treatment by the minuscule method. Table~\ref{credit}
gives a summary of which enumerations are given in which paper.

\begin{table}
\begin{center}
\begin{tabular}{|r|c||c|c|c|} \hline
\multicolumn{2}{|c||}{ }& \multicolumn{3}{l|}{Directly treated by...} \\
\multicolumn{1}{|c}{No. \cite{Stanley:symmetries}} & Case &
\multicolumn{1}{c}{This paper} &
\multicolumn{1}{c}{Stembridge \cite{Stembridge:minuscule}} &
Proctor \cite{Proctor:bruhat} \\ \hline
 1 & P    & Yes        & Yes                         & Yes \\
 2 & S    & Yes        & Yes                         & Yes \\
 3 & CS   & No         & No                          & No  \\
 4 & TS   & No         & No                          & No  \\
 5 & SC   & Yes        & Yes                         & No  \\
 6 & TC   & Yes        & No                          & No  \\
 7 & SSC  & No         & Yes                         & No  \\
 8 & CSTC & No         & No                          & No  \\
 9 & CSSC & No         & No                          & No  \\
10 & TSSC & No         & No                          & No  \\ \hline
\end{tabular}
\end{center}
\caption{\label{credit} Proctor's minuscule method and its extensions}
\end{table}

\section{Conventions}

\begin{figure}[htb]
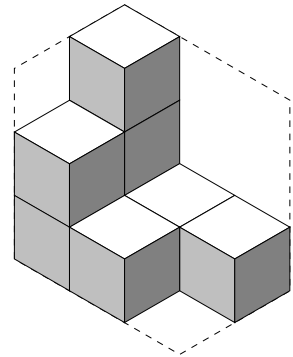

\pspicture(-2,0)(3,5.5)
\pspolygon[linestyle=dashed](0,5.5)(-1.732,4.5)(-1.732,1.5)
    (.866,0)(2.598,1)(2.598,4)
\rput(-.866,5){\lzef} \rput(-.866,5){\lzsef} \rput(.866,5){\lzswf}
\rput(.866,4){\lzswf} \rput(-1.732,3.5){\lzef} \rput(-1.732,3.5){\lzsef}
\rput(0,3.5){\lzswf} \rput(-1.732,2.5){\lzsef} \rput(0,2.5){\lzef}
\rput(-.866,2){\lzef} \rput(-.866,2){\lzsef} \rput(.866,2){\lzef}
\rput(.866,2){\lzsef} \rput(.866,2){\lzswf} \rput(2.598,2){\lzswf}
\endpspicture
\caption{A plane partition in a box}
\label{f:ppartition}
\end{figure}

A plane partition in an $a \times b \times c$ box is a collection of
unit cubes in the rectangular solid $[0,a] \times [0,b] \times [0,c]$
in three dimensions which is stable under gravitational attraction
towards the origin; \fig{f:ppartition} shows an example of a plane
partition.

The two simplest symmetry operations on the set of plane partitions are
transposition $\tau$, which consists of switching two coordinate axes; and
rotation $\rho$, which consists of cyclically permuting all three axes.  There
is also a third operation, complementation $\kappa$, which consists of taking
all cubes which are in the box but not in the diagram of the partition and
reversing all three coordinates of the box.

\begin{figure}[htb]
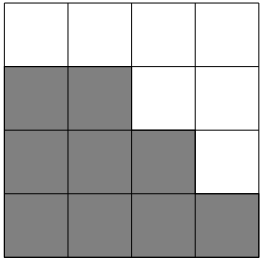

\pspicture(0,0)(4,4)
\pspolygon[fillstyle=solid,fillcolor=gray]
    (0,0)(0,3)(2,3)(2,2)(3,2)(3,1)(4,1)(4,0)
\psline(0,0)(0,4) \psline(0,0)(4,0) \psline(1,0)(1,4) \psline(0,1)(4,1)
\psline(2,0)(2,4) \psline(0,2)(4,2) \psline(3,0)(3,4) \psline(0,3)(4,3)
\psline(4,0)(4,4) \psline(0,4)(4,4)
\endpspicture
\caption{A partition in a rectangle \label{f:partition}}
\end{figure}

A partition in an $a \times b$ rectangle is the analogue of a plane partition
in two dimensions, as shown in \fig{f:partition}. The operations of
complementation and transposition act on partitions.  The set of partitions has
a partial ordering given by set-theoretic inclusion, and it will be convenient
to think of a plane partition as an ordered chain of $c$ partitions.  (Adjacent
terms may be equal in such a chain.)  Transposition of each partition in such a
chain is equivalent to transposition of the plane partition, and the same is
true of complementation provided the chain is also reversed.

A partition $P$ in an $a \times b$ rectangle can also be viewed as a a binary
sequence, i.e., a function $b_P$ from $\{1,\ldots,a+b\}$ to $\{0,1\}$ whose
values sum to $a$.  Specifically, the boundary of $P$ is a path of $a+b$ steps
from the upper left corner to the lower right corner, and we let $b_P(n) = 1$
if the $n$'th step goes down and $b_P(n) = 0$ if it goes to the right. 
Transposing or complementing a partitions has a simple effect on its binary
sequence:
\begin{eqnarray*}
b_{\kappa(P)}(n) & = & b_{P}(a+b+1-n) \\
b_{\kappa\tau(P)}(n) & = & 1-b_P(n).
\end{eqnarray*}

Finally, if $T$ is a plane partition (or a partition), $|T|$ is
the number of cubes (squares) in $T$, and the $q$-weight of $T$
is defined to be $q^{|T|}$.

\section{Arbitrary plane partitions \label{arbitrary}}

This section is a review of the following result:

\begin{theorem} The set of plane partitions in an $a \times b \times c$ box can
be regarded as a basis of an irreducible representation of $\SL(a+b)$.
\end{theorem}

This result is essentially due to Young, Schur, and Weyl, but it was  first
noted in the context of plane partitions by Stanley and Proctor, and Proctor
generalized it to symmetric plane partitions and $\SO(2a+1)$.

Consider a vector space $V$ over $\C$ (some other choices of the ground field
work also) of dimension $a+b$ and with basis $x_1,\ldots,x_{a+b}$.  The vector
space $\Lambda^a V$ of alternating forms over $V$ of degree $a$ has a basis
indexed by partitions in an $a \times b$ box.  Namely, the basis
element corresponding to $P$ is a wedge product formed from the binary sequence
of $P$:
$$x_P = \bigwedge_{n=1}^{a+b} x_n^{b_P(n)}.$$ 

Consider the vector space $S^c(\Lambda^a V)$ of symmetric forms over
$\Lambda^a V$ of degree $c$, i.e., homogeneous polynomials of degree $c$ in
the $x_P$'s, treated as independent variables.  If a plane partition $T$ is
viewed as a chain of partitions $P_1 \ge P_2 \ge \ldots \ge P_c$, then we can
define an element $p_T$ of $S^c(\Lambda^a V)$ by
$$p_T = \prod_n x_{P_n}.$$
These elements are not a basis of $S^c(\Lambda^a V)$, but merely
a linearly independent set.

We view $V$, $\Lambda^a V$, and $S^c(\Lambda^a V)$ as linear representations
of $\SL(a+b) = \SL(V)$. The representation $S^c(\Lambda^a V)$ has a unique
highest weight vector with weight $c\lambda_a$.  Here $\lambda_n$ is the
$n$th fundamental weight, the weight corresponding to the $n$th vertex from
the end of the Dynkin diagram of $\SL(a+b)$, which is a chain of length
$a+b-1$.  Therefore, the irreducible representation $V(c\lambda_a)$ appears as a
subrepresentation with multiplicity 1, and there is an equivariant projection
$$\pi:S^c(\Lambda^a V) \to V(c\lambda_a).$$
(Those not familiar with the representation theory of simple Lie algebras
should refer to Serre \cite{Serre:semisimple} or Bourbaki \cite{Bourbaki:lie}.)

We now invoke the following fundamental result:

\begin{theorem} The vectors $x_T = \pi(p_T)$, ranging
over all plane partitions $T$ in an $a \times b \times c$ box, are a
basis for the irreducible representation $V(c \lambda_a)$ of $\SL(a+b)$.
\end{theorem}

The proof follows from the theory of semistandard tableaux, originally due to
Young, as applied by Weyl to Lie groups \cite{JK:book}.

Using this result, the Weyl dimension formula gives us the
number of plane partitions in a box.  In the notation of Propp, the
number is:
$$N(a,b,c) = \frac{H(a+b+c)H(a)H(b)H(c)}{H(a+b)H(a+c)H(b+c)}$$
Here $H(n) = 1!2!3!\ldots (n-1)!$ is the hyperfactorial function.
More generally, the Weyl $q$-dimension formula gives the trace of
the element
$$D_q = q^{(1-a-b)/2} \left(\begin{array}{ccccc}
1 & 0      & 0   &        & 0 \\
0 & q      & 0   & \cdots & 0 \\
0 & 0      & q^2 &        & 0 \\
  & \vdots &     & \ddots & \vdots \\
0 & 0      & 0   & \cdots & q^{a+b-1} \end{array}\right)$$
in $\SL(a+b)$ in its action on $V(c \lambda_a)$ (or any other
irreducible representation), which is important to us because of the
following lemma:

\begin{lemma} The the action of $D_q \in \SL(a+b)$ on $V(c\lambda_a)$
is diagonal in the basis $\{x_T\}$, and the eigenvalue of
$x_T$ is $q^{|T|-abc/2}$.  \label{baby}
\end{lemma}
\begin{proof}
We will work our way up from $x_n$ for an integer $n$ to $x_P$ for a partition
$P$ to $p_T$ and then to $x_T$.  The element $D_q$ is, by definition,
diagonal in its action on $V$, and the eigenvalue of $x_n$ is $q^{n-(a+b+1)/2}$.

Since $x_P$ is a (wedge) product of $x_n$'s, it is also
an eigenvector of $D_q$, and its eigenvalue is the product of the
eigenvalues.  If $P_0$ is the empty partition, then 
$$x_{P_0} = x_1 \wedge x_2 \wedge ... \ldots x_a,$$
and the eigenvalue is therefore
$q^{a(a+1)/2-a(a+b+1)/2} = q^{-ab/2}$.  Consider the effect of an
elementary move on $P$ which consists of adding a single square to it.
This is equivalent to replacing some factor $x_n$ in $x_P$ by
$x_{n+1}$, which results in an extra factor of $q$.  Therefore,
the eigenvalue is $q^{|P|-ab/2}$ in the general case.

Since $p_T$ is a (symmetric) product of $c$ different $x_P$'s,
it is also an eigenvector of $D_q$. Its eigenvalue is the product
of the eigenvalues of its factors, namely $q^{|T|-abc/2}$.

Finally, $D_q$ commutes with the projection $\pi$ (because $\pi$ is
equivariant under the action of all of $\SL(a+b)$).  Thus, if $p_T$ is an
eigenvector with a certain eigenvalue, $x_T = \pi(p_T)$ must be an
eigenvector with the same eigenvalue.
\end{proof}

The logic of the proof of Lemma~\ref{baby} will be repeated several
more times in the rest of the paper with successively shorter
explanations.

The $q$-dimension formula yields the trace of $D_q$, which by the
lemma is $q^{-abc/2}$ times the total $q$-weight of all plane partitions.
The formula for the $q$-weight is:
$$N(a,b,c)_q = \frac{H(a+b+c)_qH(a)_qH(b)_qH(c)_q}
    {H(a+b)_qH(a+c)_qH(b+c)_q},$$
where 
$$H(n)_q = (1)_q!(2)_q!\ldots(n-1)_q!$$
is the $q$-hyperfactorial function, 
$$(n)_q! = (1)_q(2)_q\ldots(n)_q$$
is the $q$-factorial function, and
$$(n)_q = 1+q+q^2+\ldots+q^{n-1}$$
is a $q$-integer.

\section{Self-complementary plane partitions}

\begin{theorem}[Stembridge, K] If the plane partitions in an $a \times b \times
c$ box are a basis for a vector space, the operation of complementation is
linearly conjugate to the operation of negating partitions with an odd number of
cubes.
\end{theorem}

We consider the effect of the group element of the form:
$$K = i^{a+b-1}\left(\begin{array}{ccccc}
0 & 0 & 0 & 0 & 1 \\
0 & 0 & 0 & 1 & 0 \\
0 & 0 & 1 & 0 & 0 \\
0 & 1 & 0 & 0 & 0 \\
1 & 0 & 0 & 0 & 0 \end{array} \right)$$
in $\SL(a+b)$ on the $x_P$'s and the $x_T$'s.  By definition, it sends
$x_n$ to $i^{a+b-1}x_{a+b+1-n}$, and its effect on a partition is
to reverse its binary sequence, which is the complementation operation.
Specifically:
$$K(x_P) = i^{a(a+b-1)}(-1)^{a(a-1)/2}x_{\kappa(P)} = i^{ab}x_{\kappa(P)}.$$
The factor of $(-1)^{a(a-1)/2}$ comes from a reversal of the factors in
$x_P$, which is an antisymmetric product.  Since $p_T$ is a symmetric
product of $x_P$'s, the effect of $K$ is:
$$K(p_T) = i^{abc} p_{\kappa(T)}$$
Finally, since $K$ commutes with $\pi$, it has the same effect on
the $x_T$'s.  Therefore the trace of $K$ in $V(c \lambda_a)$ is, up to
a factor of $i^{abc}$, the number of plane partitions fixed by $\kappa$.
$V(b\lambda_a)$.

But $K$ is conjugate to $D_{-1}$ and therefore has the same trace.  (Note that
the notation $D_q$ is an abuse because the matrix is really a function of
$q^{1/2}$.  To obtain $D_{-1}$, we set $q^{1/2} = i$.) In conclusion,

\begin{corollary}[Stanley]
$$N_\kappa(a,b,c) = N(a,b,c)_{-1}$$
\end{corollary}
This equation is easily equivalent to the formula given by Stanley
\cite{Stanley:symmetries}.

\section{Symmetric plane partitions}

\begin{theorem} If the plane partitions in an $a \times a \times c$ box are
interpreted as a basis for an irreducible representation of $\SL(2a)$ as in
section \ref{arbitrary}, transposition is induced
by the Hodge star operator relative to a symplectic inner product.
\end{theorem}
\begin{proof}

We proceed with the notation of the section \ref{arbitrary}, taking $a = b$. 
For any non-degenerate bilinear form $B$ on $V$, the group $\SL(V)$ possesses an
outer automorphism $\sigma_B:A \mapsto A^{-T}$ (inverse transpose), where the
transpose $A \mapsto A^T$ is taken relative to the form $B$.  (Usually, $B$ is
either symmetric or antisymmetric, but it need not be either in the definition
of $\sigma_B$.) Explicitly, $B$ induces an isomorphism $\alpha_B:V \to V^*$,
where $V^*$  is the dual vector space, given by 
$$\alpha_B(v) = B(\cdot,v).$$
Then $\sigma_B(A)$ is the composition 
$$\alpha_B^{-1}{A^*}^{-1} \alpha_B,$$
where $A^*$ is the adjoint of $A$. Even more explicitly, if $A$ and $B$
are given by matrices $M_A$ and $M_B$,
$$M_{\sigma_B(A)} = M_B^{-1} M_A^{-T} M_B,$$
where $M^{-T}$ is the usual inverse transpose of the matrix $M$. We consider
the group $\Z/2 \ltimes_{\sigma_B} \SL(V)$, which is the same group for all
choices of $B$ because $\sigma_B$ differs from $\sigma_{B'}$ by an inner
automorphism of $\SL(V)$.

The automorphism $\sigma_B$ can be given a compatible linear action on
$\Lambda^a V$, called the Hodge star operator, as follows:  First, we pick a
specific volume form on $V$, i. e. an isomorphism between $\Lambda^{2a} V$
and $\C$, relative to which $B$ has determinant 1.  This identification induces
a bilinear form on $\Lambda^a V$ given by the wedge product.  There is a
second bilinear form induced by $B$ (and denoted by the same letter),
because $\Lambda^a V \subseteq V^{\tensor a}$ and $V^{\tensor a}$ has a bilinear
form which is just a multilinear extension of $B$ on $V$.  More explicitly, if
$v_1,\ldots,v_a$ and $w_1,\ldots, w_a$ are two sequences of vectors,
\begin{equation}
B(v_1 \wedge \ldots \wedge v_a,w_1 \wedge \ldots \wedge w_a) = 
\sum_{\sigma\in S_a} (-1)^{\sigma} \prod_{n=1}^a B(v_n,w_{\sigma(n)}) 
\label{expand}
\end{equation}
Both bilinear operations can be interpreted as functions from
$\Lambda^a V$ to its dual, and the Hodge star operator is either one composed
with the inverse of the other.  More explicitly, the Hodge star
is defined by the equation
$$B(\nu,* \omega) = \nu \wedge \omega$$
for $\nu,\omega \in \Lambda^a V$.
If it is assigned to $\sigma_B$,
$\Lambda^a V$ becomes a representation of $\Z/2 \ltimes_{\sigma_B} \SL(V)$.

With respect to the given basis for $V$, we define volume in the standard way,
and we define $B$ by a matrix of the form:
$$\left(\begin{array}{cccccc}
 0 & 0 & 0 & 0 & 0 & 1 \\
 0 & 0 & 0 & 0 &-1 & 0 \\
 0 & 0 & 0 & 1 & 0 & 0 \\
 0 & 0 &-1 & 0 & 0 & 0 \\
 0 & 1 & 0 & 0 & 0 & 0 \\
-1 & 0 & 0 & 0 & 0 & 0
\end{array}\right) .$$
In general, $B$ is defined by the equations
\begin{equation}
B(x_n,x_k) = \left\{
\begin{array}{ll} (-1)^n & k = 2a+1-n \\ 0 & \mbox{otherwise} \end{array}
\right. .
\label{bilin} \end{equation}
With this choice of $B$,
$$B(x_P,x_{P'}) = 0$$
unless $P' = \kappa(P)$, and
$$B(x_P,x_{\kappa(P)}) = (-1)^{|P|},$$
by multilinearity.
The sign in this equation merits some explanation.  A factor of
$(-1)^{|P|+a(a-1)/2}$ comes from multiplying $a$ factors of equation
\ref{bilin}, while another factor of $(-1)^{a(a-1)/2}$ appears because
in the sole non-zero term in equation \ref{expand}, $\sigma$ is the
reversing permutation.

On the other hand, 
$$x_P \wedge x_{P'} = 0$$
unless $b_P(n) = 1-b_{P'}(n)$, because otherwise the wedge product has
repeated linear factors.  This happens precisely when $P' = \kappa\tau(P)$.
In this case:
$$x_P \wedge x_{\kappa\tau(P)} = (-1)^{|P|}x_1 \wedge \ldots \wedge x_{2a}.$$ 
The sign comes from permuting the wedge factors.  Composing the wedge
product (as a scalar-valued bilinear function) with $B$, we obtain
$$ * x_P = x_{\tau(P)}, $$
which establishes that the Hodge star operation as the linear
extension of transposition of partitions.  We can extend the Hodge star
multilinearly
to the polynomial space $S^c(\Lambda^a V)$ over $a$-forms, and it immediately
follows that
\begin{equation}
* p_T = p_{\tau(T)}. \label{polystar}
\end{equation}

We would like to infer the conclusion of the theorem,
$$* x_T = x_{\tau(T)},$$
from equation \ref{polystar}.  To do so, we need to know that the Hodge
star operator commutes with the projection $\pi:S^c(\Lambda^a V) \to
V(c \lambda_a)$.  In general, given a representation $R$ of a
semidirect product $\Z/2 \ltimes_\sigma G$, with $G$ reductive,
either $\sigma$, acting on $G$, fixes the character of a $G$-irreducible
summand $V$ of $R$, in which case $\sigma$, acting on $R$, leaves
$V$ invariant; or $\sigma$ does not fix the character of $V$, in
which case $\sigma(V)$ is a disjoint $G$-irreducible summand of $R$.
In our case, $V(c \lambda_a)$ appears with multiplicity 1, and 
$\pi$ is characterized as the projection whose image is $V(c \lambda_a)$
and whose kernel is the direct sum of all other $SL(V)$-irreducible summands.
With the given choice of $B$, the automorphism $\sigma_B$ is induced by reversing
the Dynkin diagram of $SL(V)$, which fixes its middle fundamental
weight $\lambda_a$ and therefore also fixes $c \lambda_a$ and the
entire character of $V(c \lambda_a)$.  Therefore both the kernel
and image of $\pi$ are invariant under the Hodge star action of $\sigma_B$,
which implies that $\pi$ and $*$ commute.
\end{proof}

As stated in the proof, $\sigma_B$ is a Dynkin diagram automorphism. The
character theory of semi-direct products arising from Dynkin diagram
automorphisms is described by Neil Chriss \cite{Chriss:thesis}, who explained
to the author that although this theory is known to several representation
theorists, it may not have been previously published. The group $\Z/2
\ltimes_{\sigma_B} \SL(2a)$ has two components.  The character of a
representation on the identity component is just the usual character of
$\SL(2a)$.  The character on the $\sigma_B$ component, when non-zero, equals
the character of an associated representation of the dual Lie group, in this
case $\SO(2a+1)$, to the subgroup fixed by the outer automorphism, in this case
$\Sp(2a)$.  The representation associated to $V_{\SL(2a)}(c \lambda_a)$ is the
projective representation $V_{\SO(2a+1)}(c \lambda_a)$, where $\lambda_a$ is
now the weight corresponding to the short root of $B_a$, the root system of
$\SO(2a+1)$.  In particular, the trace of $\sigma_B$ is the dimension of
$V_{\SO(2a+1)}(c \lambda_a)$, as given by the Weyl dimension formula, and the
trace of $\sigma_BD_q$ is the $q$-dimension, as given by the Weyl $q$-dimension
formula.  From the construction of these linear operators, their traces are
$N_\tau(a,a,c)$, the number of $\tau$-invariant plane partitions, and their
total $q$-weight, $N_\tau(a,a,c)_q$, times $q^{-a^2c/2}$, respectively.

By analogy with the results described in section \ref{arbitrary}, Proctor
\cite{Proctor:bruhat} also realizes $\tau$-invariant plane partitions as a basis
for the same projective representation $V_{\SO(2a+1)}(c \lambda_a)$, and the
$q$-dimension again equals the trace of the linear transformation which
multiplies each plane partition by its $q$-weight times $q^{-a^2c/2}$.  Since
the result of this section arrives at the same Lie group character, if not an
explicit representation, it provides another proof which is similar to but not
the same as Proctor's proof of the product formulas for $N_\tau(a,a,c)$ and
$N_\tau(a,a,c)_q$. Converting the Weyl formulas by an extension of Propp's
notation, they are
$$ N_\tau(a,a,c)_q = \frac{H_2(2a+b)_{q^2} H(a)_{q^2} (2a+b-1)_q!!}
    {H(a+b)_{q^2} H_2(b+1)_{q^2} (b-1)_q!!(2a-1)_q!!} $$
Here $$H_2(n) = (n-2)_q!(n-4)_q!(n-6)_q!\ldots$$ is a staggered
$q$-hyperfactorial, and $$n_q!! = (n)_q(n-2)_q(n-4)_q\ldots$$ is
a staggered $q$-factorial.
(The formula was first proved by George Andrews.)

\section{Transpose-complement plane partitions}

In this section, we combine the results of the previous two sections to enumerate
$\kappa\tau$-invariant plane partitions.  In the representation $V(c\lambda_a)$
of $\Z/2 \ltimes_{\sigma_B} \SL(2a)$, which has plane partitions in an $a \times
a \times c$ box as a basis, the group element $\sigma_B$ acts as the linear
extension of $\tau$, while the group element $K$ acts as $i^{a^2c}$ times the
linear extension of $\kappa$.  Therefore their product acts as $i^{a^2c}$ times
$\kappa\tau$.  If $\sigma_BK$ were conjugate to  $\sigma_BD_{-1}$, it would
imply by the results of the last two sections that
$$N_{\tau\kappa}(a,a,c) = N_\tau(a,a,c)_{-1},$$
which easily implies the usual formula due to Proctor and mentioned by Stanley 
\cite{Stanley:symmetries}.  It is easy to check that
a group element with a matrix of the form
$$\frac1{\sqrt{2}} \begin{pmatrix}
 1 & 0 &  0 & 1 \\
 0 & 1 & -1 & 0 \\
 0 & 1 &  1 & 0 \\
-1 & 0 &  0 & 1 \end{pmatrix}$$
conjugates $K$ to $D_{-1}$ and commutes with $\sigma_B$.  (The general
expression for this element is $\frac1{\sqrt{2}}(I+D_{-1}K^{-1})$, or its
matrix is $\frac1{\sqrt{2}}(I+M_B)$.) Therefore $\sigma_BK$ and
$\sigma_BD_{-1}$ are also conjugate.  The conclusion is
\begin{theorem}
$$N_{\kappa\tau}(a,a,c) = N_{\tau}(a,a,c)_{-1}$$
\end{theorem}

\section{Symmetric, self-complementary plane partitions}

In this paper, we have constructed a linear representation
$V_{\SL(a+b)}(c\lambda_a)$ with a basis indexed by plane partitions.  We have 
recognized the linear extension of complementation $\kappa$ as part of the
group action on this representation.  When $a=b$, we have also recognized the
linear extension of transposition $\tau$ as part of the group action of a
semidirect product of $\SL(a+b)$. The traces of these linear transformations
tell us the number of elements fixed by $\kappa$ or $\tau$ or $\kappa\tau$. 
The author briefly thought that this construction would lead to an enumeration
of plane partitions fixed by both $\kappa$ and $\tau$.  this is a naive hope. 
Consider a permutation representation of the Klein group generated by $\tau$
and $\kappa$ with one free orbit and two fixed points, and consider a second
permutation representation with three orbits with two points, with stabilizers
$\tau$, $\kappa$, and $\kappa\tau$.  These two representations are different,
but the linear representations they induce are isomorphic. 

Thus, in the absence of extra structure, the linear representation
$V_{\SL(2a)}(c\lambda_a)$ gives no information about $\kappa,\tau$-invariant
plane partitions.  However, Stembridge enumerates $\kappa,\tau$-invariant plane
partitions by expressing $\kappa$ as an involution on the representation
$V_{\SO(2a+1)}(c\lambda_a)$. The author conjectures that latter representation
can be realized as a natural vector subspace of the former one.

\vspace{10ex}


\begin{thebibliography}{1}

\bibitem{Bourbaki:lie}
Nicolas Bourbaki, \emph{Groupes et alg\'ebres de {Lie}}, vol. 4--6, Hermann,
  Paris, 1968.

\bibitem{Chriss:thesis}
Neil Chriss, \emph{Representations of {Hecke} algebras arising from unramified
  groups}, Ph.D. thesis, University of Chicago, 1993.

\bibitem{JK:book}
Gordon James and Adalbert Kerber, \emph{The representation theory of the
  symmetric group}, Addison-Wesley, Reading, Mass., 1981.

\bibitem{Proctor:bruhat}
Robert~A. Proctor, \emph{Bruhat lattices, plane partition generating functions,
  and minuscule representations}, Europ. J. Combin. \textbf{5} (1984), no.~4,
  331--350.

\bibitem{Serre:semisimple}
Jean-Pierre Serre, \emph{Complex semismiple {Lie} algebras}, Springer-Verlag,
  New York, 1987.

\bibitem{Stanley:symmetries}
Richard~P. Stanley, \emph{Symmetries of plane partitions}, J. Combin. Theory
  Ser. A \textbf{43} (1986), no.~1, 103--113.

\bibitem{Stembridge:minuscule}
John~R. Stembridge, \emph{On minuscule representations, plane partitions, and
  involutions in complex {Lie} groups}, Duke Math. J. \textbf{73} (1994),
  no.~2, 469--490.

\bibitem{Stembridge:q-1}
\bysame, \emph{Some hidden relations involving the ten symmetry classes of
  plane partitions}, J. Combin. Theory Ser. A \textbf{68} (1994), 372--409.

\end{thebibliography}

\providecommand{\bysame}{\leavevmode\hbox to3em{\hrulefill}\thinspace}

\end{document}